# Existence and uniqueness of solutions of linear sparse matrix equations via a fixed point theorem

Xiaorong Liu*

[a]*Department of Mathematics, University of Colorado Boulder, Boulder, CO 80302, USA*

**ABSTRACT.** In this paper we prove several generalizations and applications of the Banach fixed point theorem for the complete metric spaces. This theorem is used to prove the existence and uniqueness of solutions of the linear sparse matrix problem considered.

**Keywords:** banach fixed point theorem, contraction mapping, sparse matrix.

## 1. Introduction and Preliminaries

The Banach fixed point theorem was introduced by Stefan Banach in 1922 [1, 2]. The theorem is an important tool in the theory of metric spaces in mathematics.

Let $(X, d)$ be a metric space, and $d = |\cdot|$ be the Euclidean distance. Then a map $T: X \to X$ is called a contraction mapping on $X$, if there exists $q \in [0,1)$ such that

$$d(T(x), T(y)) \leq q d(x, y)$$

holds for all $x, y$ in $X$ [1].

**Definition 1** ([3]) A mapping $T: X \to X$ is said to have a fixed point, if there exists $x \in X$ such that $T(x) = x$.

**Definition 2** ([3, 4]) Let $(X, d)$ be a non-empty complete metric space with a contraction mapping $T: X \to X$, then

(i) $T$ admits a unique fixed point $x^*$ in $X$;

(ii) $T^n(x) = x^*$ for every $x \in X$;

(iii) $d(T^n(x), x^*) \leq \frac{q^n}{1-q} d(T(x), x)$ for a real number $0 < q < 1$.

**Proof.** Let $x \in X$ be arbitrary, now we define a sequence.

$$x_0 = x, x_1 = T(x_0), x_2 = T(x_1), \ldots, x_n = T(x_{n-1}), \tag{1}$$

Since $T$ is a contraction mapping, we have

$$d(x_{n+1}, x_n) = d(T(x_n), T(x_{n-1})) \leq q d(x_n, x_{n-1}),$$

for $0 < q < 1$, according to the definition of contraction mapping and the sequence (1).

After repeating this argument for n times, we get

$$d(x_{n+1}, x_n) \leq qd(x_n, x_{n-1}) \leq q^2 d(x_{n-1}, x_{n-2}) \leq \cdots \leq q^n d(x_1, x_0).$$

Now, apply the triangle inequality [5] into this formula. For any positive integer $p$, we have

$$\begin{aligned} d(x_{n+p}, x_n) &\leq d(x_{n+p}, x_{n+p-1}) + d(x_{n+p-1}, x_{n+p-2}) + \cdots + d(x_{n+1}, x_n) \\ &\leq q^{n+p-1} d(x_1, x_0) + q^{n+p-2} d(x_1, x_0) + \cdots + q^n d(x_1, x_0) \\ &= \sum_{i=1}^{p} q^{n+i-1} d(x_1, x_0), \end{aligned}$$

and

$$\sum_{i=1}^{p} q^{n+i-1} = q^n \sum_{i=0}^{p-1} q^i = q^n \frac{(1-q^i)}{1-q} = \frac{q^n}{1-q}.$$

Since $q < 1$, $q^i \to 0$ as $i \to \infty$. Hence

$$d(x_{n+p}, x_n) \leq \frac{q^n}{1-q} d(x_1, x_0) \tag{2}$$

When $p$ gets large enough.

If $n \to \infty$, $q^n \to 0$. Let $\varepsilon > 0$ be arbitrary, from (2) it follows that

$$d(x_{n+p}, x_n) \leq \frac{\varepsilon}{1-q} d(x_1, x_0) = \varepsilon.$$

When $n$ is large enough, this proves that the sequence (1) is Cauchy. Thus, there is a point $x^* \in X$ such that

$$\lim_{n \to \infty} x_n = x^*.$$

Now, take the limit on the both sides of the equation $x_n = T(x_{n-1})$, we obtain $\lim_{n \to \infty} x_n = \lim_{n \to \infty} T(x_{n-1})$. Since $T$ is continuous, we can take the limit inside. Then we see that

$$x^* = \lim_{n \to \infty} x_n = \lim_{n \to \infty} T(x_{n-1}) = T\left(\lim_{n \to \infty} x_{n-1}\right) = T(x^*), \tag{3}$$

so $x^*$ is a fixed point in $X$.

Suppose there exists another fixed point $y^*$ in $X$ such that

$$T(y^*) = y^*.$$

According to the definition of contraction mapping, it follows

$$d(x^*, y^*) = d(T(x^*), T(y^*)) \leq qd(x^*, y^*),$$

and hence $(1-q)d(x^*, y^*) \leq 0$. Remembering $q \in [0,1)$, it follows that $0 \leq 1 - q \leq 1$. So we have $d(x^*, y^*) = 0$, i.e. $y^* = x^*$. Thus, we proved (i).

Since x was arbitrary, from (3) it follows that

$$x^* = \lim_{n\to\infty} x_n = \lim_{n\to\infty} T(x_{n-1}) = \lim_{n\to\infty} T(T(x_{n-2})) = \lim_{n\to\infty} T \ldots T(x) = T^n(x),$$

then we have proved (ii).

From (2), we take $p \to \infty$, then

$$d(x^*, T^n(x)) = d(x_{n+p}, x_n) \leq \frac{q^n}{1-q} d(x_1, x_0) = \frac{q^n}{1-q} d(T(x), x),$$

so (iii) is proved. And the proof of the theorem is complete.

Recently, it has been shown that the well-known Banach Fixed Point Theorem is useful in many applications [5]. In this paper, we give applications to integral equations, matrix equations, and to sparse matrix equations.

## 2. The main results

**Theorem 3** Consider the metric space $(X, d)$ where $X = R$. Let $f: [a, b] \to [a, b]$ be a differential equation that has a Lipschitz constant less than 1. Then there exists a unique $x^* \in [a, b]$ such that

$$f(x^*) = x^*.$$

**Proof.** For any $x, y \in [a, b]$ where $x < y$, we have by the mean value theorem [6] that there exists a point $c$ in an open interval $(x, y)$ at which

$$f(y) - f(x) = f'(c)(y - x).$$

After taking the absolute value both sides of the equation, we have

$$|f(y) - f(x)| = |f'(c)||y - x| \leq k|y - x|$$

for some positive $k \leq 1$, which implies that $f$ is a contraction mapping with a bounded factor $k$. By Banach fixed point theorem, there exists a unique $x^* \in [a, b]$ such that $f(x^*) = x^*$.

**Example 4.** ([7]) Define a function $f: [0,1] \to R$ which satisfies $\int_0^1 \frac{\sin(f(t)-y)}{2} dy = f(t)$ for $t \in [0,1]$. Let $M = C([0,1])$ and the mapping $A: M \to M$ be a function with

$$A(f)(t) = \int_0^1 \frac{\sin(f(t) - y)}{2} dy,$$

then there exists a unique function $f^*$ which satisfies $A(f^*) = f^*$.

**Proof.** Suppose $k$ is Lipschitz continuous to the first variable with Lipschitz constant 1 and $k(x, y) = \sin(x - y)$,

$$A(f)(t) = \int_0^1 \frac{k(f(t), y)}{2} dy.$$

Then, we have

$$|A(f)(t) - A(g)(t)| \le \frac{1}{2}\int_0^1 |k(f(t),y) - k(g(t),y)|\,dy = \frac{1}{2}\int_0^1 |k(f(t)-g(t), y-y)|\,dy$$
$$= \frac{1}{2}\int_0^1 |\sin(f(t)-g(t))|\,dy \le \frac{1}{2}|f(t)-g(t)|$$

for all $t \in [0,1]$.

So $|A(f) - A(g)| \le \frac{1}{2}|f-g|$ which shows that the mapping $A$ is a contraction with a Lipschitz constant less than 1. Based on theorem 3, there exists a unique function $f^*$ which satisfies $A(f^*) = f^*$.

**Theorem 5** Let $a, b \in R$ with $a < b$ and $(x,y) \to M(x,y)$ be a measurable function on $\{a \le x \le b, a \le y \le b\}$. Suppose $\int_a^b \int_a^b |M(x,y)|^2 dxdy < \infty$ and the function $f, g$ are square-integrable function, then the integral equation

$$f(x) = g(x) + u\int_a^b M(x,y)f(y)dy$$

has a unique solution for $|u| < |M(x,y)|^{-1}$.

In mathematics, a square-integrable function [8, 9] is a real measurable function for which the integral of the square of the absolute value is finite. Thus, if

$$\int_{-\infty}^{+\infty} |\varphi(x)|^2 dx < \infty,$$

then $\varphi$ is square integrable on the real line $(-\infty, +\infty)$ with the following property

$$\|\varphi\| = \left(\int_S |\varphi|^2 d\mu\right)^{1/2}$$

for $(S, \mu)$ which is a measure space.

**Proof.** Since $\int_a^b \int_a^b |M(x,y)|^2 dxdy < \infty$, the function $M(x,y)$ is a square-integrable function. Let a function $F$ defined by

$$F(x) = g(x) + u\int_a^b M(x,y)f(y)dy$$

We need to show that $\varphi(x) = \int_a^b M(x,y)f(y)dy$ is also a square-integrable function. By Holder's inequality [10], we have

$$\int_a^b |\varphi(x)|^2 dx = \int_a^b \left|\int_a^b M(x,y)f(y)dy\right|^2 dx \le \int_a^b \left(\int_a^b |M(x,y)|\,|f(y)|\,dy\right)^2 dx$$
$$\le \int_a^b \left(\int_a^b |M(x,y)|^2\,dy\right) \left(\int_a^b |f(y)|^2\,dy\right) dx.$$

From Fubini's theorem [11], we further get

$$\int_a^b \left( \int_a^b |M(x,y)|^2 \, dy \right) \left( \int_a^b |f(y)|^2 \, dy \right) dx = \left( \int_a^b \int_a^b |M(x,y)|^2 \, dydx \right) \left( \int_a^b |f(y)|^2 \, dy \right) < \infty.$$

Thus, we obtain $\int_a^b |\varphi(x)|^2 dx < \infty$. That is, $\varphi(x) = \int_a^b M(x,y)f(y)dy$ is a square-integrable function.

Define a mapping $T(f) = F$ where $T$ is also a square-integrable function. By the property of square-integrable function and Holder's inequality, we have

$$d(T(f_1), T(f_2)) = \left( \int_a^b |F_1 - F_2|^2 dx \right)^{1/2}$$

$$= |u| \left( \int_a^b \left| \int_a^b M(x,y)f_1(y)dy - \int_a^b M(x,y)f_2(y)dy \right|^2 dx \right)^{1/2}$$

$$= |u| \left( \int_a^b \left| \int_a^b M(x,y)(f_1(y) - f_2(y))dy \right|^2 dx \right)^{1/2}$$

$$\leq |u| \left( \int_a^b \left( \int_a^b |M(x,y)|^2 dy \right) \left( \int_a^b |f_1(y) - f_2(y)|^2 dy \right) dx \right)^{\frac{1}{2}}$$

$$= |u| \left( \int_a^b \left( \int_a^b |M(x,y)|^2 dy \right) dx \right)^{\frac{1}{2}} \left( \int_a^b |f_1(y) - f_2(y)|^2 dy \right)^{\frac{1}{2}}$$

$$= |u| \left( \int_a^b \int_a^b |M(x,y)|^2 dy \, dx \right)^{\frac{1}{2}} d(f_1, f_2).$$

Since $|u| < |M(x,y)|^{-1}$, $T$ is a contraction mapping. By Definition 2, there exists a unique $f^*$ such that $(f^*) = f^*$.

**Theorem 6** Consider the matrix problem $Ax = b$, where

$$A = \begin{pmatrix} a_{11} & \cdots & a_{1n} \\ \vdots & \ddots & \vdots \\ a_{n1} & \cdots & a_{nn} \end{pmatrix}, \quad x = (x_1, \cdots x_n)^T, \quad b = (b_1, \cdots b_n)^T.$$

Then, there exists a unique solution if

$$|I - A| \leq q$$

for $0 < q < 1$.

**Proof.** We can re-write this system of equations as

$$x_1 = (1 - a_{11})x_1 - a_{12}x_2 - \cdots - a_{1n}x_n + b_1;$$

$$x_2 = -a_{21}x_1 + (1 - a_{22})x_2 - \cdots - a_{2n}x_n + b_2;$$

$$\vdots$$

$$x_n = -a_{n1}x_1 - a_{n2}x_2 - \cdots + (1 - a_{nn})x_n + b_n.$$

For $1 \leq i, j \leq n$, set $\alpha_{ij} = (I - A) = \begin{bmatrix} 1 - a_{11} & -a_{12} & \cdots & -a_{1n} \\ -a_{21} & 1 - a_{22} & \cdots & -a_{2n} \\ \vdots & & \ddots & \vdots \\ -a_{n1} & -a_{n2} & \cdots & 1 - a_{nn} \end{bmatrix}$. So the above system of equations can be written as

$$x_i = \sum_{j=1}^{n} \alpha_{ij} x_j + b_i.$$

We can find the matrix problem $Ax = b$ is equivalent to the matrix problem $x = x - Ax + b$. Define a function $T: R^n \to R^n$ which satisfies $T(x) = x - Ax + b$. Thus, for $x_1, x_2 \in R^n$,

$$T(x_1) - T(x_2) = x_1 - x_2 - (Ax_1 - Ax_2) = (I - A)(x_1 - x_2)$$

According to the definition of matrix, we have

$$d(T(x_1), T(x_2)) = \sqrt[2]{\sum_{j=1}^{n} \left(\alpha_{ij}(x_{1j} - x_{2j})\right)^2} \leq \sqrt[2]{\sum_{j=1}^{n} (\alpha_{ij})^2 \left(x_{1j} - x_{2j}\right)^2} \leq \left|\sum_{j=1}^{n} \alpha_{ij}\right| \sqrt[2]{\sum_{j=1}^{n} \left(x_{1j} - x_{2j}\right)^2}$$

$$= \left|\sum_{j=1}^{n} \alpha_{ij}\right| d(x_1, x_2) = |I - A| d(x_1, x_2) \leq q d(x_1, x_2).$$

which shows that $T$ is a contraction mapping. The proof of the theorem is complete using Definition 2.

**Example 7.** Consider the matrix problem $Ax = b$, where

$$A = \begin{pmatrix} 3/2 & 1 & 0 \\ 0 & 1 & 1 \\ 0 & 1 & 2/3 \end{pmatrix}, \quad x = (x_1 \quad x_2 \quad x_3)^T, \quad b = (2 \quad 3 \quad 1)^T,$$

it follows that

$$|I - A| = \left|\begin{pmatrix} -1/2 & -1 & 0 \\ 0 & 0 & -1 \\ 0 & -1 & 1/3 \end{pmatrix}\right| = (-1) \times \left|\begin{pmatrix} -1/2 & -1 & 0 \\ 0 & -1 & 1/3 \\ 0 & 0 & -1 \end{pmatrix}\right| = 1/2 \leq q$$

for $0 < q < 1$. Then, there exists a unique solution $x = (10/3 \quad -3 \quad 6)^T$.

For this specific case, if we determine whether the solution is unique by traditionally computing the determinant of A, it is easy to confirm that the computation complexity is greater than using Theorem 6 above. This illustrates that Theorem 6 can decrease the computation burden of determining the uniqueness of the solution to matrix equations.

**Lemma 8.** Assume that $Ax = b$, and $A$ is a sparse matrix such that

$$A = \begin{pmatrix} a_{11} & \cdots 0 \cdots & 0 \\ \vdots & & \vdots \\ 0 & \ddots & 0 \\ \vdots & & \vdots \\ 0 & \cdots 0 \cdots & a_{nn} \end{pmatrix}, \quad x = (x_1, \cdots x_n)^T, \quad b = (b_1, \cdots b_n)^T.$$

Then there exists a unique solution if $|I - A| \leq q$ for $0 < q < 1$.

As we know, the sparse matrix is a matrix in which most of the elements are zeros. Thus, this lemma satisfies all conditions of Theorem 6. So the proof of Lemma 8 is completed according to Theorem 6. In general, the standard operations and algorithms of the sparse matrix equations are very slow since the processing is wasted on the zeros. Based on Lemma 8, we can determine the existence and uniqueness of the solutions of sparse matrix problems by only one step.

**Example 9.** Suppose $x = (1,1,1,0,1,1)^T$ is one solution of matrix problem $Ax = b$, where

$$A = \begin{pmatrix} 11 & 22 & 0 & 0 & 0 & 0 \\ 0 & 0 & 44 & 0 & 0 & 0 \\ 0 & 0 & 2 & 66 & 0 & 0 \\ 0 & 0 & 0 & 0 & 77 & 0 \\ 0 & 0 & 0 & 0 & 88 & 99 \\ 0 & 0 & 0 & 0 & 0 & 1/1111 \end{pmatrix}, x = (x_1, \cdots x_6)^T, b = (33 \quad 44 \quad 2 \quad 77 \quad 187 \quad 1/1111)^T.$$

Determine the all solutions of this equations.

Notice that $|I - A| < 871/1111$. Thus, all the conditions guarantee that $x = (1,1,1,0,1,1)^T$ is the only solution.